\documentclass[12pt,twoside]{article} 
\usepackage[cp1251]{inputenc}
\usepackage{latexsym,amsfonts,amssymb,amsmath}
\usepackage{euscript}
\usepackage[T2A]{fontenc}
\usepackage{amsthm}
\usepackage[english]{babel}
\usepackage{mathrsfs}

\begin{document}


\newcommand{\ang}[1]{\langle{#1}\rangle}

\thispagestyle{empty}

\noindent {UDC 519.41/47}

\setcounter{page}{1}

\vspace{1mm}\noindent{\bf N.S.~Chernikov} (In-te Math. of NAS of
Ukraine, Kyiv, Ukraine)

\vspace{1mm}\noindent{\bf GROUPS WITH THE MINIMAL CONDITIONS FOR\\
SUBGROUPS AND FOR NONABELIAN SUBGROUPS}

\noindent chern@imath.kiev.ua

\vspace{0.5cm}\noindent{\small For some very wide classes
$\mathfrak{D}$ and $\mathfrak{B}\subset\mathfrak{D}$ of groups, the
author proves that an arbitrary (nonabelian) group $G\in
\mathfrak{D}$ (respectively $G\in \mathfrak{B}$) satisfies the
minimal condition for (nonabelian) subgroups iff it is Chernikov.}

\thispagestyle{empty}

\vspace{8mm}

\noindent Recall that a group $G$ is called Shunkov, if for any its
finite subgroup $K$ every subgroup of the quotient group
$N_{G}(K)/K$, generated by two its conjugated elements of prime
order, is finite. Recall that a group is called locally graded, if
any its finitely generated  subgroup $\neq 1$ contains a subgroup of
finite index $\neq 1$ (S.N.Chernikov). The class of all periodic
Shunkov groups is wide and contains, for instance, all binary finite
and 2-groups. The class of all locally graded groups is very wide
and contains, for instance, all abelian, locally finite, residually
finite groups. It is easy to see that this class is local and any
group having a series with locally graded factors is locally graded,
and any residually (locally graded) group is locally graded. At the
same time, the class of all locally graded groups contains all
groups having a series with locally finite factors, all $RN$- (and
so all groups of all Kurosh-Chernikov classes), locally solvable,
locally hyperabelian, radical in the sense of B.I.Plotkin,
residually solvable groups.

Below, as usual, $\pi(G)$ is the set of all primes $p$ for which the
group $G$ has a $p$-element $\neq 1$.

Let $\mathfrak{A}$ be the class of all groups $G$ for which the
following conditions are fulfilled:
\begin{enumerate}
    \item [(i)] If $G$ is not torsion-free, then for any $p\in
    \pi(G)$ and $p$-element $g\neq 1$ of $G$ and $h\in
    C_{G}(g^{p})$, $\ang{g,g^{h}}$ possesses a subgroup of finite
    index $\neq 1$ or $\ang{g,h}$ is not periodic.
    \item [(ii)]If $G$ is not periodic, then for any element $g$ of
    infinite order of $G$ and $h\in G$, $\ang{g,g^{h}}$ possesses a subgroup of finite
    index $\neq 1$.
\end{enumerate}

The class $\mathfrak{A}$ is very wide and contains, for instance,
all periodic Shunkov groups and all locally graded groups.

Let $\mathfrak{C}$  be the class of all groups $G$ for which (i),
with "periodic $\neq 1$" instead of "not torsion-free" and deleted
"or $\ang{g,h}$ is not periodic", is fulfilled (and (ii) is not
necessary fulfilled). Then $\mathfrak{A}$ is contained in
$\mathfrak{C}$, $\mathfrak{C}$ contains all nonperiodic groups and
the class of all periodic $\mathfrak{A}$-groups is just the class of
all periodic $\mathfrak{C}$-groups.

Let $\mathfrak{B}$ (resp. $\mathfrak{D}$) be the minimal local class
of groups containing $\mathfrak{A}$ (resp. $\mathfrak{C}$) such that
any group possessing a series with $\mathfrak{B}$ (resp.
$\mathfrak{D}$)-factors belongs to $\mathfrak{B}$ (resp.
$\mathfrak{D}$). Put $\mathfrak{B_{0}}=\mathfrak{A}$ (resp.
$\mathfrak{D}_{0}=\mathfrak{C}$) and for ordinals $\beta>0$ by
induction: if for some ordinal $\alpha$, $\beta=\alpha+1$, then
$\mathfrak{B_{\beta}}$ (resp. $\mathfrak{D_{\beta}}$) be the class
of all groups which have a local system of subgroups possessing a
series with $\mathfrak{B_{\alpha}}$ (resp.
$\mathfrak{D_{\alpha}}$)-factors, and if there is no such $\alpha$,
then $\mathfrak{B_{\beta}}={\mathop {\bigcup}
\limits_{\alpha<\beta}} \mathfrak{B_{\alpha}}$ (resp.
$\mathfrak{D_{\beta}}={\mathop {\bigcup} \limits_{\alpha<\beta}}
\mathfrak{D_{\alpha}}$). For some ordinal $\gamma$,
$\mathfrak{B}=\mathfrak{B}_{\gamma}$ (resp.
$\mathfrak{D}=\mathfrak{D}_{\gamma}$). It is easy to show by
induction that all $\mathfrak{B_{\beta}}$ and, at the same time,
$\mathfrak{B}$ are closed with respect to subgroups.

The known Shunkov's \cite{1} and S.N.Chernikov's \cite{2} Theorems
establish that a nonabelian group satisfying the minimal condition
for nonabelian subgroups is Chernikov, if it is locally finite or
has a series with finite factors resp. The next Theorem contains
them.

\vskip2mm \textbf{Theorem.} \textit{A (nonabelian) group $G\in
\mathfrak{D}$ (resp. $G\in \mathfrak{B}$) satisfies the minimal
condition for (nonabelian) subgroups iff it is Chernikov.}

Note that Ol'shanskii's nonabelian groups, in which all proper
subgroups are finite (see \cite{3}), satisfy the minimal condition
for subgroups and are non-Chernikov. Thus, in Theorem, the condition
"$G\in \mathfrak{D}$" is essential. Note that Ol'shanskii's
nonabelian torsion-free groups, in which all proper subgroups are
cyclic (see \cite{3}), satisfy the minimal condition for nonabelian
subgroups, are Shunkov and non-Chernikov. In particular, in Theorem,
the condition ''$G\in \mathfrak{B}$'' is essential.

Below $\min$ and $\min - \overline{ab}$ are the minimal conditions
for subgroups and nonabelian subgroups resp. Other notations are
standard. (Remark that a group with $\min$ is periodic and an
abelian group with $\min$ is Chernikov.)

\vspace{3mm}{\itshape\bfseries Proof. }  {\it Sufficiency } is
obvious.

\hspace{0.2cm}{\it Necessity. } Let $G$ be non-Chernikov.

(a) {\it Reduction to the case when $G$ satisfies $\min -
\overline{ab}$ and also $G\in \mathfrak{A}$.} Let $\zeta$ be minimal
among all ordinals $\alpha$, for which $\mathfrak{B}_{\alpha}$
(resp. $\mathfrak{D}_{\alpha}$) contains a non-(Chernikov or
abelian) group (resp. a non-Chernikov group) with $\min -
\overline{ab}$ (resp. $\min$). We may assume that
$G\in\mathfrak{B_{\zeta}}$ (resp. $G\in\mathfrak{D_{\zeta}}$).
Suppose $\zeta>0$. Clearly, for some ordinal $\xi$, $\zeta=\xi+1$.
So $G$ possesses a local system of subgroups having a series with
$\mathfrak{B}_{\xi}$ (resp. $\mathfrak{D}_{\xi}$)-factors. Every
factor is, obviously, Chernikov or abelian (resp. Chernikov). So $G$
is locally graded. Thus $G\in \mathfrak{B}_{0}$ (resp. $G\in
\mathfrak{D}_{0}$), which is a contradiction. So $\zeta=0$ and $G\in
\mathfrak{A}$ (resp. $G\in \mathfrak{C}$). In the case of $\min$,
$G\in \mathfrak{C}$ and $G$ is, obviously, periodic nonabelian with
$\min - \overline{ab}$. In particular, $G\in \mathfrak{A}$. Taking
this into account, we may consider later on only the case of $G\in
\mathfrak{A}$ with $\min - \overline{ab}$.

Since $G$ satisfies $\min - \overline{ab}$, it contains a
non-(Chernikov or abelian) subgroup $L$ such that any its proper
subgroup is Chernikov or abelian. We may assume that $G=L$. Then for
every normal subgroup $N$ of $G$, any proper subgroup of $G/N$ is
Chernikov or abelian.

(b) {\it Show that a subgroup $H$ of $G$ is Chernikov or abelian, if
it has a subgroup $K$ of finite index $\neq 1$ or if it is almost
solvable.} Since $K\neq G$, $K$ and so $H$ are almost abelian. If
$H$ is almost solvable, then in view of Corollary to Theorem 1
\cite{2}, Corollary 2 \cite{2} and Lemmas 1,2 \cite{2}, it is
Chernikov or abelian.

(c) {\it Show that $G$ is periodic.} Let $G$ have elements $g$ of
infinite order. Since $G\in \mathfrak{A}$, every $\langle
g,g^{h}\rangle$ has a proper subgroup of finite index and so is
abelian (see (b)). Then for any $u\in G$, $\ang{g^{h}:h\in
G}\ang{u}$ is non-Chernikov solvable and so is abelian (see (b)).
Thus $g\in Z(G)$. Let $v\in G$ and $|\ang{v}|<\infty$. Then $vg$ is
of infinite order. So $v=(vg)g^{-1}\in Z(G)$. Thus $G$ is abelian,
which is a contradiction.

(d) {\it Show that $G/Z(G)$ is Shunkov.} Let $K/Z(G)$ be a finite
subgroup of $G/Z(G)$. In view of Kalu\v{z}nin's Theorem (see
\cite{4}), $C_{G}(K/Z(G))/C_{G}(K)$ is abelian. So
$N_{G}(K)/C_{G}(K)$ is almost abelian.

If $K/Z(G)\neq 1$, then $C_{G}(K)\neq G$ and so $C_{G}(K)$ is almost
abelian. Thus $N_{G}(K)$ is almost solvable. So $N_{G}(K)$ is
Chernikov or abelian (see (b)). Clearly,
$N_{G}(K)/Z(G)=N_{G/Z(G)}(K/Z(G))$. Consequently, the quotient group
$N_{G/Z(G)}(K/Z(G))/(K/Z(G))$ is Chernikov or abelian. Therefore any
two its elements of prime order generate a finite subgroup.

Let $K/Z(G)=1$ and $R/Z(G)$ be a subgroup of $G/Z(G)$ generated by
two its conjugated element of some prime order $p$. Obviously,
because of $Z(G)$ is periodic (see (c)), for some $p$-element$g\in
G$ such that $g^{p}\in Z(G)$ and some $h\in G$,
$R=\ang{g,g^{h}}Z(G)$. Clearly, $R/Z(G)$ is isomorphic to a quotient
group of the group $\ang{g,g^{h}}/\ang{g^{p}}$. Since $G\in
\mathfrak{A}$, $\ang{g,g^{h}}$ has a subgroup of finite index $\neq
1$. So, with regard to (b), $\ang{g,g^{h}}$ is finite. Therefore
$R/Z(G)$ is finite.

(e) {\it Show that $G/Z(G)$ has an abelian non-Chernikov nonnormal
maximal subgroup $A/Z(G)$ such that
\begin{equation}\label{1}
    A/Z(G)\cap (A/Z(G))^{g}=1,\quad g\in (G/Z(G))\backslash(A/Z(G)).
\end{equation}}If all proper subgroups of $G/Z(G)$ are Chernikov, it satisfies
$\min$. So because of $G/Z(G)$ is Shunkov (see (d)), by
Suchkova-Shunkov Theorem \cite{5} it is Chernikov. So $G$ is almost
solvable, which is a contradiction (see (b)). So some maximal
abelian subgroup $A/Z(G)$ of $G/Z(G)$ is non-Chernikov. An arbitrary
proper subgroup $H\supseteq A$ of $G$ is non-Chernikov. So it is
abelian. Therefore $H/Z(G)=A/Z(G)$ and $H=A$. Thus $A$ is an abelian
maximal subgroup of $G$. If $A$ is normal in $G$, then $|G:A|$ is
prime and $G$ is solvable, which is a contradiction (see (b)).
Consequently, for any $g\in G\backslash A$, $G=\ang{A,A^{g}}$. Then
$A\cap A^{g}\subseteq Z(G)$. But, clearly, $Z(G)\subseteq A,A^{g}$.
Thus $A\cap A^{g}= Z(G)$. Therefore (1) is valid.

(f) {\it Show that $A/Z(G)$ has some element $a$ of odd prime
order.} Suppose that this is not the case. Since $A/Z(G)$ is
periodic (see (c)) and neither cyclic nor quasicyclic, it has some
elements $b$ and $c\neq b$ of order 2. Let $h\in (G/Z(G))\backslash
(A/Z(G))$. Then
$\ang{b,c^{h}}=\ang{bc^{h}}\leftthreetimes\ang{b}=\ang{bc^{h}}\leftthreetimes\ang{c^{h}}$
and $|\ang{bc^{h}}|<\infty$ (see (c)). If $|\ang{bc^{h}}|$ is odd,
then for some $s\in \ang{bc^{h}}$, $b=c^{hs}$. Since $A/Z(G)$ is
abelian and $b,c\in A/Z(G)$, $hs\notin A/Z(G)$. But $b\in A/Z(G)\cap
(A/Z(G))^{hs}$, which is a contradiction (see (1)). So
$\ang{bc^{h}}$ contains some element $w$ of order 2. But then $w\in
C_{G/Z(G)}(b)\cap C_{G/Z(G)}(c^{h})=A/Z(G)\cap(A/Z(G))^{h}$, which
is a contradiction.

(g) {\it Final contradiction.} Let $a$ be from (f). Since $G/Z(G)$
is Shunkov (see (d)), for any $h\in G/Z(G)$,
$|\ang{a,a^{h}}|<\infty$. So with regard to (1) by Sozutov-Shunkov
Theorem \cite{6}, for some normal subgroup $N/Z(G)$ of $G/Z(G)$,
$G/Z(G)=(A/Z(G))(N/Z(G))$ and $A/Z(G)\cap N/Z(G)=1$. Since $N\neq G$
and $G/N$ is abelian, $G$ is almost solvable, which is a
contradiction (see (b)).

\vspace{3mm} The following proposition is contained in Theorem.

\vskip2mm \textbf{Proposition 1.} \textit{Let $G$ be a (nonabelian)
group. Assume that $G$ is a periodic Shunkov or locally graded or 2-
group. Then $G$ satisfies $\min$ (resp. $\min-\overline{ab}$) iff it
is Chernikov.}

In view of Mal'cev Theorem (see Theorem 4.2 \cite{8}), a linear
group over a field is locally residually finite. Further, for a
commutative and associative ring $R$ with 1 and any finitely
generated unital module $M$ over $R$, $\mathrm{Aut}_{R}(M)$ is
hyperabelian-by-residually (linear over fields) (Theorem 13.5
\cite{8}). Consequently, $\mathrm{Aut}_{R}(M)$ is locally graded.
Hence follows that any ${\rm\mathbf{GL}}_{n}(R)$ is locally graded.
Therefore, in virtue of Theorem, the following proposition is valid.

\vskip2mm \textbf{Proposition 2.} \textit{A (nonabelian) group
$G\subseteq \mathrm{Aut}_{R}(M)$ or $G\subseteq
{\rm\mathbf{GL}}_{n}(R)$ satisfies $\min$ (resp.
$\min-\overline{ab}$) iff it is Chernikov.}

\vspace{0.5cm}Finally, let $\mathfrak{E}$ be the class of all groups
$G$ for which the following conditions are fulfilled:
\begin{enumerate}\itemsep0.1cm
  \item [(i)] If $G$ is not torsion-free, then for any $p\in \pi(G)$
  and $p$-element $g\neq 1$ of $G$ and $h\in C_{G}(g^{p})$,
  $\ang{g,g^{h}}$ possesses a $\mathfrak{B}$-homomorphic image $\neq
  1$ or $\ang{g,h}$ is not periodic.
  \item [(ii)] If $G$ is not periodic, then for any element $g$ of
    infinite order of $G$ and $h\in G$, $\ang{g,g^{h}}$ possesses a
    $\mathfrak{B}$-homomorphic image $\neq 1$.
\end{enumerate}

Let $\mathfrak{F}$ be the class of all groups $G$ for which the
following condition is fulfilled:

If $G$ is periodic $\neq 1$, then for any $p\in \pi(G)$ and
$p$-element $g\neq 1$ of $G$ and $h\in C_{G}(g^{p})$,
$\ang{g,g^{h}}$ possesses a $\mathfrak{D}$-homomorphic image $\neq
  1$.

\vskip2mm \textbf{Proposition 3.} \textit{A (nonabelian) group $G\in
\mathfrak{F}$ (resp. $G\in \mathfrak{E}$) satisfies $\min$ (resp.
$\min-\overline{ab}$) iff it is Chernikov.}

\vspace{3mm}{\itshape\bfseries Proof. } {\it Necessity.} A
corresponding homomorphic image of $\ang{g,g^{h}}$ is generated by
two elements. In the case when it is not abelian, by Theorem, it is
finite. In the case when it is abelian, it has a subgroup of finite
index $\neq 1$. Consequently, $\ang{g,g^{h}}$ possesses a subgroup
of finite index $\neq 1$. So $G\in \mathfrak{C}$ (resp. $G\in
\mathfrak{A}$). Therefore by Theorem, $G$ is Chernikov.

{\it Sufficiency } is obvious.

\end{document}